\documentclass[10pt]{article}
\usepackage{amssymb}
\usepackage{amsfonts}
\usepackage{amsmath}

\newtheorem{theorem}{Theorem}

\newtheorem{lemma}[theorem]{Lemma}

\newtheorem{proposition}[theorem]{Proposition}

\newenvironment{proof}[1][Proof]{\noindent\textbf{#1.} }{\ \rule{0.5em}{0.5em}}

\begin{document}

\title{Increasing hazard rate of mixtures for natural exponential families}
\author{ Shaul K. Bar-Lev$\thanks{%
Department of Statistics, University of Haifa, Haifa 31905, Israel
(barlev@stat.haifa.ac.il)}$ \ and G\'{e}rard Letac$\thanks{%
Laboratoire de Statistique et Probabilit\'es, Universit\'e Paul Sabatier,
31062 Toulouse, France (gerard.letac@alsatis.net) }$}
\maketitle

\begin{abstract}
Hazard rates play an important role in various areas, e.g., reliability
theory, survival analysis, biostatistics, queueing theory and actuarial
studies. Mixtures of distributions are also of a great preeminence in such
areas as most populations of components are indeed heterogeneous. In this
study we present a sufficient condition for mixtures of two elements\ of the
same natural exponential family (NEF) to have an increasing hazard rate. We
then apply this condition to some classical NEF's having either quadratic,
or cubic variance functions (VF) and others as well. A particular attention
is devoted to the hyperbolic cosine NEF having a quadratic VF, the Ressel
NEF having a cubic VF and to the Kummer distributions of type 2 NEF. The
application of such a sufficient condition is quite intricate and
cumbersome, in particular when applied to the latter three NEF's. Various
lemmas and propositions are needed then to verify this condition for these
NEF's.

\end{abstract}

\section{Introduction}

Hazard rates (also called failure rates) play an important role in various
areas, e.g., reliability theory, queuing models, survival analysis and
actuarial studies. Mixtures of distributions are also of a great preeminence
in such areas as most populations of components are indeed heterogeneous. A
comprehensive list of references on the behavior of hazard rates for
mixtures of distributions can be found in the monograph by Shaked and
Shanthikumar (2007) and the references cited therein and also in Block et
al. (2003).

In this study we present \ (Proposition 2) a sufficient condition for
mixtures of two elements of the same natural exponential family (NEF) to
have an increasing hazard rate. We then apply this condition to some
classical absolutely continuous NEF's having either quadratic or cubic
variance functions (VF's) (c.f., Morris, 1982, Letac and Mora, 1990) and
others NEF's as well. A particular attention is devoted to the hyperbolic
cosine NEF having a quadratic VF, the Ressel NEF having a cubic VF and the
Kummer distributions of type 2 NEF. The application of such a sufficient
condition can be intricate, in particular when applied to the latter three
NEF's. Various lemmas and propositions are then needed to verify this
condition for such NEF's. Accordingly, we dedicate Sections 4, 5 and 6,
respectively, for these three NEF's. In Section 3 we consider the rather
easy application of the sufficient condition to three NEF's having either
quadratic or cubic VF's, namely, the normal, gamma and the inverse Gaussian
NEF's. Our sufficient condition stems from the following seminal result by
Glaser (1980):

\begin{proposition}
Suppose that the probability density $s(x)$, concentrated on the interval $%
(a,\infty )$ (with $-\infty \leq a<\infty $), is positive such that $%
-b(x)=\log s(x)$ is concave. Then the mapping $x\mapsto \log
\int_{x}^{\infty }s(t)dt$ is concave on $(a,\infty )$ and the hazard
function $h(x)=s(x)/\int_{x}^{\infty }s(t)dt$ is increasing.
\end{proposition}

(Ron Glaser observes for the one line proof that since $b^{\prime }(x)$ is
nondecreasing, one has $(1/h)^{\prime }(x)=\int_{x}^{\infty
}e^{b(x)-b(t)}(b^{\prime }(t)-b^{\prime }(x))dt\leq 0$). If the probability
density is $s=e^{-b}$ the fact that $b$ is convex is by Proposition 1 a
sufficient condition for having $h$ increasing but \textit{not} a necessary
one: see the remark in Section 2 introducing the Glaser set as well the
Jorgensen and the Karlin sets, or consider the density $s_2=e^{-b_2}$ in Section 5 for which $h$ in increasing and $b_2$ is not convex.

Our sufficient condition for a mixture of two members in the same NEF to
have an increasing hazard rate is as follows: Suppose that the NEF is
written as 
\begin{equation*}
\{e^{-\lambda x-k(\lambda )-b(x)}\mathbf{1}_{(a,\infty )}(x)dx,\text{ }%
\lambda \in \Lambda \},
\end{equation*}%
where $\Lambda $ is an interval and $-\infty \leq a.$ Suppose also that $%
b^{\prime \prime }(x)\geq 0$ for all $x>a$ and denote $T(x)=1/\sqrt{%
b^{\prime \prime }(x)}$. We show in Proposition 2 that if there exists $c>0$
and $d\in \mathbb{R}$ such that the inequality $cT(x)\leq \cosh (cx+d)$
holds for all $x>a$, we then can find pairs $\lambda _{1}$ and $\lambda
_{2}=\lambda _{1}+2c$ in $\Lambda $ and a mixing coefficient $p\in (0,1)$
such that the mixture density 
\begin{equation*}
\left( pe^{-\lambda _{1}x-k(\lambda _{1})}+(1-p)e^{-\lambda _{2}x-k(\lambda
_{2})}\right) e^{-b(x)}
\end{equation*}%
has an increasing hazard rate. This simple condition relies on the fact that
the mixture is employed with two elements of the same NEF. However, this two
element mixture result is apparently not extendable to a more multi-element
mixture situation.

\section{A sufficient condition for mixtures of members of the same NEF to
have an increasing hazard rate}

Consider an absolutely continuous NEF concentrated on $(a,\infty )$ with $%
-\infty \leq a<\infty $, and generated by a locally integrable function $s$
on $(a,\infty )$ such that 
\begin{equation}  \label{LKT}
L(\lambda )=e^{k(\lambda)}=\int_{a}^{\infty }e^{-\lambda x}s(x)dx,
\end{equation}
the Laplace transform (LT) of $s(x)$, exists on a nonempty open interval $%
\Lambda $. The corresponding NEF is then given by the set of probability
densities on $(a,\infty )$ of the form%
\begin{equation}
\left\{ \exp \left\{ -\lambda x-k(\lambda) \right\}s(x) dx,\text{ }\lambda
\in \Lambda \right\} .  \label{Exp}
\end{equation}%
Let $\nu (d\lambda )$ be a probability on $\Lambda $ and suppose that the
function on $(a,\infty )$ defined by 
\begin{equation}
R(x)=\int_{\Lambda }e^{-\lambda x}\frac{\nu (d\lambda )}{L(\lambda )}
\label{R}
\end{equation}%
exists. Thus, $s(x)R(x)dx$ is a probability density on $(a,\infty )$ and it
is a mixture of the elements of the NEF. This probability density has the
hazard rate 
\begin{equation}
h_{\nu }(x)=\frac{s(x)R(x)}{\int_{x}^{\infty }s(t)R(t)dt}.  \label{hnu}
\end{equation}%
Proposition 1 shows that $h_{\nu }(x)$ is increasing if $x\mapsto \log
(s(x)R(x))$ is concave, or, equivalently, if $s(x)>0$ for all $x>a$, if $%
s^{\prime \prime }(x)$ exists and if on $(a,\infty )$ one has 
\begin{equation}
\frac{s^{\prime \prime }(x)s(x)-\left( s^{\prime }(x)\right) {}^{2}}{s^{2}(x)%
}+\frac{R^{\prime \prime }(x)R(x)-\left( R^{\prime }{}(x)\right) ^{2}}{%
R^{2}(x)}\leq 0.  \label{H}
\end{equation}
We now have the following proposition when $\nu $ is a mixture of two Dirac
measures.

\begin{proposition}
Consider the special case of the hazard rate $h_{\nu }$ in (\ref{hnu}) with%
\begin{equation}
\nu =p\delta _{\lambda _{1}}+(1-p)\delta _{\lambda _{2}},  \label{nu}
\end{equation}%
where $p\in (0,1)$ and $\lambda _{1}<\lambda _{2}$ with $\lambda _{1}$ and $%
\lambda _{2}$ in $\Lambda .$ Assume that on $(a,\infty )$, $s(x)>0$, $%
-b(x)=\log s(x)$ is concave and that $s^{\prime \prime }(x)$ exists and
define 
\begin{equation}
T(x)=1/\sqrt{b^{\prime\prime}(x)} \label{T}
\end{equation}%
\begin{equation}
p_{1}=pe^{-k(\lambda _{1})},\ p_{2}=(1-p)e^{-k(\lambda _{2})},\ c=\frac{%
\lambda _{2}-\lambda _{1}}{2}\text{ and }\ d=\log \sqrt{p_{1}/p_{2}}.
\label{EP}
\end{equation}%
Then the hazard rate (\ref{hnu}) with $\nu $ as in (\ref{nu}) is increasing
if for all $x>a$ 
\begin{equation}
cT(x)\leq \cosh (cx+d).  \label{P}
\end{equation}
\end{proposition}

\begin{proof}
The proof is a straightforward application of (\ref{H}). Indeed, for $\nu $
in (\ref{nu}) and $R$ defined by (\ref{R}), 
\begin{equation*}
R(x)=pe^{-k(\lambda _{1})}e^{-\lambda _{1}x}+(1-p)e^{-k(\lambda
_{2})}e^{-\lambda _{2}x}=p_{1}e^{-\lambda _{1}x}+p_{2}e^{-\lambda _{2}x},
\end{equation*}%
implying that $R^{\prime \prime }(x)R(x)-\left( R^{\prime }(x)\right)
{}^{2}=p_{1}p_{2}(\lambda _{2}-\lambda _{1})^{2}e^{-(\lambda _{1}+\lambda
_{2})x}$. Accordingly, the inequality (\ref{H}) becomes for this particular
case 
\begin{equation}
p_{1}p_{2}(\lambda _{2}-\lambda _{1})^{2}e^{-(\lambda _{1}+\lambda
_{2})x}\leq (p_{1}e^{-\lambda _{1}x}+p_{2}e^{-\lambda _{2}x})^{2}\left( -%
\frac{s^{\prime }(x)}{s(x)}\right) ^{\prime },\text{ }x>a.  \label{HG}
\end{equation}%
Since $\log s(x)$ is concave, $(-s^{\prime }(x)/s(x))^{\prime }\geq 0$ on $%
(a,\infty )$, so that $T(x)$ in (\ref{T}) is well defined. Thus, with the
notations (\ref{EP}), the inequality (\ref{HG}) is equivalent to 
\begin{equation*}
2cT(x)=T(x)(\lambda _{2}-\lambda _{1})\leq \sqrt{\frac{p_{1}}{p_{2}}}e^{-%
\frac{\lambda _{1}-\lambda _{2}}{2}x}+\sqrt{\frac{p_{2}}{p_{1}}}e^{-\frac{%
\lambda _{2}-\lambda _{1}}{2}x}=2\cosh (cx+d).
\end{equation*}%
which is (\ref{P}).
\end{proof}

In the next sections, we are going to consider a number of absolutely
continuous NEF's on the real line, generated by a density $s$, and check for
each of them whether $T$ exists or not, that is whether $s$ is log concave
or not. When $T$ exists we will have to discover which $(c,d)$ with $c>0$
are such that (\ref{P}) holds for all $x.$ As we shall see, for some NEF's
such that $s$ is log concave it may occur that (\ref{P}) does not hold for
any $(c,d).$

The system (\ref{EP}) of equalities links the three parameters $(\lambda
_{1},\lambda _{2},p)$ with the two parameters $(c,d).$ Suppose that we are
given a pair $(c,d)$ satisfying (\ref{P}), we therefore may choose
arbitrarily the mean $\lambda =\frac{\lambda _{1}+\lambda _{2}}{2}$ in $%
\Lambda $ such that $\lambda _{1}=\lambda -c$ and $\lambda _{2}=\lambda +c$
are in $\Lambda .$ Having done this choice of $\lambda $, the value of of
the mixing coefficient $p$ in (\ref{EP}) can be determined exactly as 
\begin{equation}
p=\frac{e^{d}L(\lambda -c)}{e^{d}L(\lambda -c)+e^{-d}L(\lambda +c)}\text{,}
\label{PP}
\end{equation}%
where $L$ is the Laplace transform (LT) of the generating density $s$
(recall that $s$ is not necessarily a probability). We note, however, that
the LT is not always expressible in terms of simple functions but rather in
terms of transcendental or implicit functions, in which case a numerical
search is then needed to find the $(c,d)$ interval on which the appropriate
mixture density possesses an increasing hazard rate. As this paper is rather
theoretical, we do not intend to pursue such a numerical search.

\vspace{4mm}\noindent \textbf{Remarks on the Jorgensen, Karlin and
Glaser sets.} Given a density $s$ on $(a,\infty )$ with LT (\ref{LKT}) such
that $\Lambda $ is not empty, the set $J(s)$ of $\alpha \geq 0$ such that $%
L^{\alpha }$ is still a LT of a positive measure $\mu _{\alpha }$ is called
the Jorgensen set of $s$ (see for instance Letac and Mora (1990) and the
references cited therein). By definition, $J(s)$ is a closed additive
semigroup. Note that $s$ generates an NEF of infinitely divisible
distributions if and only if $J(s)=[0,\infty ).$ If not, $J(s)$ can be
complicated. For instance, a consequence of the short and elegant paper by
Ben Salah and Masmoudi (2010) is that $J(s)=[1,\infty )$ if 
\begin{equation*}
s(x)=\frac{1}{4}e^{-x}\mathbf{1}_{(0,\infty )}(x)+\frac{3}{4}e^{-x+1}\mathbf{%
1}_{(1,\infty )}(x).
\end{equation*}%
Note also that $\mu _{\alpha }$ could have a continuous singular part for some
small $\alpha \in J(s)$ although appropriate examples are rather
complicated. If $J^{\ast }(s)\subset J(s)$ is the set of $\alpha $ such that 
$\mu _{\alpha }(dx)=s_{\alpha }(x)dx$ has density we trivially have that $%
J^{\ast }(s)+J(s)\subset J^{\ast }(s)$, since the convolution of a measure
with a density with any measure has a density.

Consider now a non trivial result due to Karlin and Proschan (1960) (see
also Karlin, 1968, p. 152 and Barlow and Proshan, 1965, p. 100) which says
that if $s$ and $\ell $ are probability densities with increasing hazard
rate, then the convolution $s\ast \ell $ has the same property. Therefore
let us introduce the Karlin set $K(s)$ of $\alpha \in J^{\ast }(s)$ such
that $s_{\alpha }$ has an increasing hazard rate. The above property shows
that $K(s)$ is a closed additive subsemigroup of $J^{\ast }(s).$ For
instance if $s(x)=e^{-x}\mathbf{1}_{(0,\infty )}(x)$, it is a simple
exercise to see that $K(s)=[1,\infty ).$

Finally, consider the Glaser set $G(s)$ of $s$ which is the set of $\alpha $
in the Karlin set $K(s)$ for which the conditions of Proposition 1 are met,
that is such that if $s_{\alpha }=e^{-b_{\alpha }}$ then $b_{\alpha }$ is
convex. Although in many cases $G(s)$ coincides with $K(s)$ the Glaser set $G(s)$ is not a
semigroup. Indeed, we shall face with an example in Section 4 relating to
the Ressel distribution $s_{1}$ where $G(s)$ is a bounded interval; thus
distinct from the semi group $K(s).$ Usually the Karlin set is more difficult to find than the Glaser set. 

\section{Applications related to NEF's with quadratic or cubic VF's (normal,
gamma and inverse Gaussian NEF's)}

As already noted in the introduction, quadratic VF's include six NEF's of
which only three have densities: Normal, gamma and hyperbolic cosine (c.f.
Morris, 1982). Cubic VF's include also six NEF's of which only two have
densities: Inverse Gaussian and Ressel (c.f. Letac and Mora, 1990). The set
of our examples will include all of the five absolutely continuous NEF's
having either quadratic or cubic VF's and also another one, the Kummer
distributions of type 2 NEF. The present section deals with the normal,
gamma and inverse Gaussian NEF's and Sections 4, 5 and 6 consider the three
other ones. In what follows and whenever feasible, we provide, for each of
the examples, with their respective VF $(V,\Omega )$, where $V$ is the VF
corresponding to (\ref{Exp}) and $\Omega $ is the domain of means.\bigskip

\textbf{Example 1:} \textit{The normal NEF}

The normal NEF has a constant VF, i.e. $(V,\Omega )=(\sigma ^{2},\mathbb{R})$%
. For a fixed standard deviation $\sigma $, the generating density is 
\begin{equation*}
s(x)=\frac{1}{\sqrt{2\pi }\sigma }e^{-\frac{x^{2}}{2\sigma ^{2}}},\text{ }%
a=-\infty .
\end{equation*}%
Trivially here the Glaser set $G(s)$ is $(0,\infty).$ This leads to $%
T(x)=\sigma $ and $k(\lambda )=\frac{\sigma ^{2}\lambda ^{2}}{2}.$ The
inequality (\ref{P}) is fulfilled for any $x\in \mathbb{R}$ if and only if $%
c\sigma \leq 1$, or equivalently, if $|\lambda _{1}-\lambda _{2}|\leq
2/\sigma $, a result that was already obtained by Block \textit{et al.}
(2005).\bigskip

\textbf{Example 2:} \textit{The gamma NEF}

The gamma NEF, concentrated on $\left( 0,\infty \right) $, has a VF $V(\mu
)=\alpha^{-1}\mu ^{2}$ and $\Omega =\mathbb{R}^{+}$, where $\alpha $ and $%
\mu $ are, respectively, the shape and mean parameters. For a fixed shape
parameter $\alpha >0$, the generating measure is 
\begin{equation*}
s_{\alpha}(x)dx=\frac{x^{\alpha -1}}{\Gamma (\alpha )}\mathbf{1}%
_{(0,\infty)}(x)dx.
\end{equation*}%
We now consider three exhausted cases relating to the values of the
parameter $\alpha$: $\alpha =1$ (the exponential case), $\alpha<1$ and $%
\alpha >1$. These observations imply that the Glaser set $G(s_1)$ is $%
[1,\infty).$

\begin{enumerate}
\item $\alpha =1$. Here, $s(x)\equiv 1$ so that for any $\nu $ the function $%
\log R$ is convex and the inequality (\ref{P}) (as well as (\ref{H})) cannot
be fulfilled unless $\nu $ is concentrated on one point.

\item $\alpha <1$. Since for this case both $\log R$ and $\log s$ are
convex, the inequalities (\ref{H}) or (\ref{P}) cannot be fulfilled.

\item $\alpha >1$. Here, $T(x)=x/\sqrt{\alpha-1}$ and we have the following
proposition.
\end{enumerate}

\begin{proposition}
For $\alpha >1$, the probability density 
\begin{equation*}
f(x)=\frac{1}{\Gamma (\alpha )}x^{\alpha -1}(p\lambda _{1}^{\alpha
}e^{-\lambda _{1}x}+(1-p)\lambda _{2}^{\alpha }e^{-\lambda _{2}x}),x>0,
\end{equation*}
where $\lambda_1<\lambda_2,$ has an increasing hazard rate $%
h(x)=f(x)/\int_{x}^{\infty }f(t)dt$ if 
\begin{equation*}
\frac{\lambda _{2}}{\lambda _{1}}\leq \left( \frac{p}{1-p}\right) ^{1/\alpha
}e^{-d_{0}/2\alpha },
\end{equation*}%
where
\end{proposition}

\begin{equation*}
d_{0}=\log \frac{1+\sqrt{\alpha }}{\sqrt{\alpha -1}}-\frac{2\alpha }{(1+%
\sqrt{\alpha })}.
\end{equation*}%
\vspace{4mm}\noindent \textbf{Example.} For $\alpha =2$ the result
specializes to the following: since $d_{0}=\log 2-2$ and since $%
e^{d_{0}/2}=e/\sqrt{2}=1.92...$ we can claim that the following mixing of
two gamma densities 
\begin{equation*}
f(x)=p\lambda _{1}^{2}xe^{-\lambda _{1}x}+(1-p)\lambda _{2}^{2}xe^{-\lambda
_{2}x}
\end{equation*}%
where $\lambda _{1}<\lambda _{2}$ has an increasing hazard rate if $\left( 
\frac{\lambda _{2}}{\lambda _{1}}\right) ^{2}\leq \frac{p}{1-p}(1.92..)$.
For instance choosing $\lambda _{2}=2\lambda _{1}$ imposes a heavy weight $p$
on $\lambda _{1}$, namely $.675<p<1.$

\begin{proof}
\vspace{4mm}\noindent\ Since $c>0$ for studying the inequality (\ref{P}) we
write $t=cx+d.$ Thus (\ref{P}) becomes: for all $t>-d$, 
\begin{equation}  \label{GA}
\frac{t-d}{\sqrt{\alpha -1}}\leq \cosh t.
\end{equation}%
For a fixed $\alpha >1$, we determine the set of $d$ values such that (\ref%
{GA}) holds. Since $t\mapsto \cosh t$ is a convex function, we look for the
point $(t_{0},\cosh t_{0})$ such that the tangent to the curve $\cosh $ has
slope $1/\sqrt{\alpha -1}.$ Thus, 
\begin{equation*}
\sinh t_{0}=\frac{1}{\sqrt{\alpha -1}},\ t_{0}=\log \frac{1+\sqrt{\alpha }}{%
\sqrt{\alpha -1}}\text{ and}\ \cosh t_{0}=\frac{2\alpha }{(1+\sqrt{\alpha })%
\sqrt{\alpha -1}}.
\end{equation*}%
The equation of this tangent is $y=\left( t-d_{0}\right) /\sqrt{\alpha -1},$
where $d_{0}$ is such that this line goes through the point $(t_{0},\cosh
t_{0}).$ This implies that 
\begin{equation}
d_{0}=t_{0}-(\cosh t_{0})\sqrt{\alpha -1}=\log \frac{1+\sqrt{\alpha }}{\sqrt{%
\alpha -1}}-\frac{2\alpha }{(1+\sqrt{\alpha })}.  \label{B}
\end{equation}%
Such results show that $t/\sqrt{\alpha -1}\leq \cosh (t+d)$ for all $t>0$ if
and only if $d\geq d_{0}.$ The application of this fact is that (\ref{HG})
holds for all $x>0$ if and only if $\sqrt{p_{1}/p_{2}}\geq e^{d_{0}}$, or,
equivalently, if 
\begin{equation*}
\frac{\lambda _{1}}{\lambda _{2}}\geq \left( \frac{1-p}{p}\right) ^{1/\alpha
}e^{d_{0}/2\alpha }.
\end{equation*}%
\medskip
\end{proof}

\textbf{Example 3:} \textit{The inverse Gaussian NEF} The inverse Gaussian
NEF has a VF $V(\mu )=\alpha ^{-2}\mu ^{3}$ with $\Omega =\mathbb{R}^{+}$,
where $\alpha >0$ and $\mu $ is the mean parameter. Here $a=0$\ and for a
fixed $\alpha >0$, the corresponding NEF is generated by 
\begin{equation*}
s(x)=\frac{\alpha }{\sqrt{2\pi }}x^{-3/2}e^{-\frac{\alpha ^{2}}{2x}}.
\end{equation*}%
This implies that $(-s^{\prime }(x)/s(x))^{\prime }=(3x-2\alpha
^{2})/2x^{3}, $ which is not a positive function. Thus Proposition 2 is not
applicable.\medskip

\section{The\ hyperbolic cosine NEF}

The hyperbolic cosine NEF $\mathcal{F}_{\alpha }$ has a VF $V(\mu )=\frac{%
\mu ^{2}}{\alpha }+\alpha $ with $\Omega =\mathbb{R}$, where $\alpha >0$.
The convex support of $\mathcal{F}_{\alpha }$ is $\mathbb{R}$ (i.e., $%
a=-\infty $). The generating measure of $\mathcal{F}_{\alpha }$ is 
\begin{equation*}
\mu _{\alpha }(dx)=\frac{2^{\alpha -2}}{\pi }\left\vert \Gamma (\frac{\alpha
+ix}{2})\right\vert ^{2}\frac{dx}{\Gamma (\alpha )}.
\end{equation*}%
(See Morris 1982, for details). Moreover, its LT is defined on $\Lambda =(-%
\frac{\pi }{2},\frac{\pi }{2})$ by $L_{\alpha }(\lambda )=(\cos \lambda
)^{-\alpha }$ (since it is $1$ for $\lambda =0$, this shows the non obvious
fact that $\mu _{\alpha }$ is a probability). The hyperbolic cosine $\mu
_{1} $ distribution is more known as the hyperbolic secant (hereafter, HS
distribution or HS NEF). Various probabilistic properties of the HS
distribution have been derived, though it is rarely used in applied
statistics, probably due to its intricate structure. Although this
distribution is not used much in applications, it does, however, have two
curious features: Like the normal distribution, the density of $\mu _{1}$ is
proportional to its characteristic function; the sample mean and median are,
asymptotically, equally efficient. A probabilistic interpretation  of $\mu _{1}$ is available: consider 
 a standard complex Brownian motion $Z=X+iY $ with $Z(0)=0$  and 
the hitting time $T$ of the set $\{x+iy\ ;|y|\geq \pi /2\}.$ Then $X(T)\sim
\mu _{1}:$ to see this, consider the process $M(t)=\exp sZ(t).$
Since $z\mapsto e^{sz}$ is analytic, it is harmonic, $M$ is a martingale and 
$\mathbb{E}(M(T))=1$ gives the desired result. A newsworthy statistical
analysis and data fitness can be found in Smyth (1994) and recently in
Sibuya (2006) (a complete English version of the latter paper is available
by corresponding the author).

Denote by $s_{\alpha}=e^{-b_{\alpha}}$ the density of $\mu _{\alpha}.$ The fact that the  function $b_{\alpha}$ is convex if and only if $\alpha\geq 1$
has been proved by Shanbhag (1979). We give a  different proof in the following proposition: 
\begin{proposition} The function $b_{\alpha}$ is convex if and only if $\alpha\geq 1$ (in other terms the Glaser set $G(s_1)$ is $[1,\infty)).$ More specifically
for $\alpha>1$ we have 
\begin{equation}  \label{HC}
\frac{1}{s_{\alpha}(x)}=(\alpha-1)\int_{-\pi/2}^{\pi/2}e^{xu}(\cos
u)^{\alpha-2}\; du.
\end{equation}
and for $\alpha<1$ the function $b''_{\alpha}$ is negative in the interval $$\left(\alpha\sqrt{\frac{2+\alpha}{2-\alpha}}\ ,\ \frac{2+\alpha}{\sqrt{3}}\right).$$
\end{proposition}

\begin{proof}
Formula (\ref{HC}) is the particular case $\nu =\alpha -1>0$ and $a=ix$ of
the classical formula 
\begin{equation}\label{GR}
\int_{0}^{\pi /2}(\cos u)^{\nu -1}\cos au\;du=\frac{\pi }{2^{\nu }\nu B(%
\frac{\nu +1+a}{2},\frac{\nu +1-a}{2})}
\end{equation}%
which can be found in Gradshteyn and Ryzhik (1980), page 372, 3.631 formula
9. Now (\ref{HC}) shows that $1/s_{\alpha }$ is the Laplace transform of the
positive measure 
\begin{equation}
(\alpha -1)(\cos u)^{\alpha -2}\mathbf{1}_{(-\pi /2,\pi /2)}(u)\;du
\label{HD}
\end{equation}%
which implies that the function $b_{\alpha }=-\log s_{\alpha }$ is strictly
convex for $\alpha>1.$ 
For $\alpha=1$ we can see that $b_1=-\log s_1$ is convex by the same trick since $1/s_1=2\cosh \frac{\pi x}{2}$ is the Laplace transform of the positive measure 
$\delta_{-\pi/2}+\delta_{\pi/2}$. Note that  is the weak limit of (\ref{HD}) when $%
\alpha\rightarrow 1.$ 

Suppose now that $0<\alpha<1.$ We use the digamma function $\psi =\Gamma ^{\prime }/\Gamma $ and its derivative. If $z$ is a complex number with positive real part :%
\begin{equation}
\psi ^{\prime }(z) =\sum_{n=0}^{\infty }\frac{1}{(n+z)^{2}}\label{DG}
\end{equation}%
A easy calculation leads to 
$$b''_{\alpha}(x)=\frac{1}{4}\psi'(\frac{\alpha+ix}{2})+\frac{1}{4}\psi'(\frac{\alpha-ix}{2})=\frac{1}{4}\varphi_{\alpha/2}(\frac{x^2}{4})$$
where
 for $t>0$ we define $$\varphi_{c}(t)=\frac{c^2-t}{\left(c^2+t\right)^2}+\sum_{n=1}^{\infty }\frac{(n+c)^2-t}{\left((n+c)^2+t\right)^2}.$$
For showing that for $c<1/2$ the function $t\mapsto \varphi_{c}(t)$ is negative on some interval observe that for fixed $t>0$ the function on $(0,\infty)$ defined by $u\mapsto \frac{u-t}{(u+t)^2}$ is decreasing when $u>3t.$ As a consequence if $(1+c)^2>3t$ we can write the majorization of the sum of a series by an integral
$$\varphi_{c}(t)<\frac{c^2-t}{\left(c^2+t\right)^2}+\int_0^{\infty}\frac{(v+c)^2-t}{\left((v+c)^2+t\right)^2}dv=\frac{c^2-t}{\left(c^2+t\right)^2}+\frac{c}{c^2+t}$$ 
(here we have used $\frac{v^2-t}{(v^2+t)^2}=-\frac{d}{dv}\frac{v}{v^2+t}).$ This shows that $\varphi_{c}(t)<0$ when $\frac{c^2(1+c)}{1-c}\leq t\leq \frac{(1+c)^2}{3}.$ Since $c<1/2$ we have$\frac{(1+c)^2}{3}-\frac{c^2(1+c)}{1-c}=\frac{1+c}{3(1-c)}(1-4c^2)>0$ and this interval is not empty. Replacing $c$ by $\alpha/2$ and $t$ by $x^2/4$ we get that  $b''_{\alpha}$ is negative in the interval indicated in the statement of the proposition.

\end{proof}

For part $0<\alpha<1$ our proof of Proposition 4 is elementary. For part $\alpha\geq 1$ our  proof  is based on the formula \ref{GR} and Laplace transforms. The compact and  ingenious Shanbhag's proof rather relies on Fourier transforms through the formula  $b''_{\alpha}(x)=\int _{-\infty}^{\infty}e^{itx/2}d_{\alpha}(t)dt$
where $$ d_{\alpha}(t)=\frac{t}{2\sinh \frac{t}{2}}e^{\frac{|t|}{2}(1-\alpha)}.$$ This formula is derived from an integral formula for $b_{\alpha}$ which can be found in Zolotarev (1967) and which is obtained from the L\'evy measure of the infinitely divisible distribution of $\log X$ when $X$ is $\gamma_{\alpha/2}$ distributed.  If $\alpha\geq 1$  the function $d_{\alpha}$ is an integrable characteristic function (corresponding to a Cauchy distribution with parameter $(\alpha-1)/2$ convoluted with  the density $\pi/2(\cosh \pi x)^2)).$ The Fourier inversion formula shows that $b''_{\alpha}(x)\geq 0$ for all $x.$ If $0<\alpha<1$ we have  $d_{\alpha}(x)>1$ around zero and $d_{\alpha}$ cannot be a characteristic function. This prevents $b''_{\alpha}$ to be positive by a careful but standard reasoning using again the Fourier inversion and this concludes the Shanbhag's proof.

Since the functions $b_\alpha$ and $T=1/\sqrt{b^{\prime\prime}_{\alpha}}$
are not simple when $\alpha$ is not an integer we therefore emphasize the
analysis of the respective mixtures for the two cases: $\alpha=1$ and $%
\alpha=2$. In principle an analysis similar to the case $\alpha=2$ below
could be also performed for $\alpha=3,4,\ldots$ but the case $\alpha=2$ is
creative enough to let us think that higher cases are difficult. 

{\textbf{The hyperbolic case }$\mathbf{\alpha=1}$}\textbf{: }\noindent The
most popular member of the $\mathcal{F}_\alpha$'s is related to this case.
As mentioned in the proof of Proposition 4 the corresponding density for $\alpha=1$ is
\begin{equation*}
s_{1}(x)=\frac{1}{2\cosh \frac{\pi x}{2}}.
\end{equation*}%
 Thus the above results are applicable to this $%
\mathcal{F}_{1}$. More specifically 
\begin{equation*}
b_1^{\prime\prime}(x)= \left( -\frac{s_{1}^{\prime }(x)}{s_{1}(x)}\right)
^{\prime }=(\frac{\pi }{2})^{2}\frac{1}{\cosh ^{2}\frac{\pi x}{2}}>0,
\end{equation*}%
and thus $T(x)=\frac{2}{\pi }\cosh \frac{\pi x}{2}$. In order to study the
inequality (\ref{P}) for this particular case, we use the following lemma.

\begin{lemma}
\noindent Let $a$ and $u$ be positive numbers and $v$ be a real number.
Then, the following inequality 
\begin{equation*}
a\cosh x\leq \cosh (ux+v)
\end{equation*}%
holds for all real $x$ if and only if $a\in (0,1],$ $u\geq 1$ and $|v|\leq
v_{0}=v_{0}(a,u)$, where 
\begin{equation*}
v_{0}=u\log \,(\frac{A}{a}+\frac{uB}{a})-\log \,(A+B)
\end{equation*}%
with the notation $A=\sqrt{\frac{u^{2}-a^{2}}{u^{2}-1}}$ and $B=\sqrt{\frac{%
1-a^{2}}{u^{2}-1}}.$
\end{lemma}

\begin{proof}
\vspace{4mm} $\Rightarrow .$ Letting $x\rightarrow \infty $, we have that $%
a\cosh x\sim ae^{x}$ and $\cosh (ux+v)\sim e^{ux}$, implying that $u\geq 1.$
Letting $x=-v/u$ shows that $a\cosh (-u/v)\leq 1$ and thus $a\leq 1.$ In the
sequel we assume that $u>1$ and treat the case $u=1$ separately after. Now,
we introduce the two positive numbers $x_{0}$ and $v_{0}$ such that the two
curves $x\mapsto a\cosh x$ and $x\mapsto \cosh (ux-v_{0})$ are tangent on a
point of the abscissa $v_{0}.$ Thus, they satisfy the two equations: 
\begin{equation*}
a\cosh x=\cosh (ux_{0}-v_{0})\text{ and}\ \ a\sinh x=u\sinh (ux_{0}-v_{0}).
\end{equation*}%
Squaring these two equations and using the fact that $\cosh ^{2}t-\sinh
^{2}t=1$, we get a linear system in $\cosh ^{2}x_{0}$ and $\cosh
^{2}(ux_{0}-v_{0})$, whose solution is 
\begin{equation*}
\cosh ^{2}x_{0}=\frac{A^{2}}{a^{2}}\text{ and}\ \cosh
^{2}(ux_{0}-v_{0})=A^{2}.
\end{equation*}%
Since $t\geq 0$ and $y=\cosh t$, it follows that $t=\log (y+\sqrt{y^{2}-1}.$
Thus, $x_{0}=\log \,(\frac{A}{a}+\frac{uB}{a})$ and $v_{0}=u\log \,(\frac{A}{%
a}+\frac{uB}{a})-\log \,(A+B)$ (note that $ux_{0}-v_{0}\geq 0$). To end the
proof of $\Rightarrow $ we show that $a\cosh x\leq \cosh (ux+v)$ for all
real $x$ would imply that $|v|\leq v_{0}.$ Since the function $v\mapsto
\cosh (ux_{0}-v)-a\cosh x_{0}$ is decreasing on the interval $(-\infty
,ux_{0})$ and is zero on $v_{0}$ (which belongs to this interval), we get
that $v<v_{0}$ when $\cosh (ux_{0}-v)-a\cosh x_{0}\geq 0.$ Similarly,
because of the symmetry of $\cosh t$ we have $-v_{0}\leq v.$

We now prove $\Leftarrow .$ Assume that $a\leq 1<u$ and $|v|\leq v_{0}.$
Denoting $f(x)=\cosh (ux-v)-a\cosh x$, then since $f(x_{0})=f^{\prime
}(x_{0})=0$, the Taylor formula gives 
\begin{equation}
f(x)=\int_{x_{0}}^{x}(x-t)f^{\prime \prime }(t)dt.  \label{Z}
\end{equation}%
We use the latter formula to show that $f(x)>0$ for $x>x_{0}.$ Note that
since $f^{\prime \prime }(x)>f(x)$, then $f^{\prime \prime }(x_{0})>0$ and (%
\ref{Z}) implies that $f(x)>0$ on some interval $(x_{0},x_{1}).$ Now suppose
that there exists $x_{2}>x_{0}$ such that $f(x_{2})=0.$ Without loss of
generality we assume that $f(x)>0$ on $(x_{0},x_{2})$. Thus $f^{\prime
\prime }(x)>0$ on $(x_{0},x_{2}).$ Since by (\ref{Z}), $f(x_{2})=0$ is
impossible, we obtain that $f(x)>0$ for all $x>x_{0}.$ To prove that $f(x)>0$
for all $x<x_{0}$ is similar. \noindent

We now consider the particular case\textsc{\ }$u=1.$ The inequality $a\cosh
x\leq \cosh (x+v)$ is equivalent to $e^{2x}(e^{v}-a)\geq a-e^{-v}$. By
letting $x\rightarrow \pm \infty $, it can be easily seen that the latter
inequality holds for all $x$ if and only if $|v|\leq -\log a.$
\end{proof}

\vspace{4mm}\noindent We do not apply the full strength of this lemma for
our problem, but instead study when the inequality 
\begin{equation*}
\frac{2c}{\pi }\cosh \frac{\pi x}{2}\leq \cosh (cx+d),
\end{equation*}%
holds for all $x$.

To fit with the notations of the latter lemma, denote $t=\frac{\pi x}{2}$
which leads to $\frac{2c}{\pi }\cosh t\leq \cosh (\frac{2c}{\pi }x+d).$ Now,
the lemma implies that if this inequality holds for all $t$ then $a=\frac{2c%
}{\pi }\leq 1\leq u=\frac{2c}{\pi }.$ Thus $\frac{2c}{\pi }=1$, but since $c=%
\frac{1}{2}(\lambda _{2}-\lambda _{1})$ we must have $\lambda _{2}=\lambda
_{1}+\pi .$ However, this is impossible since the corresponding LT $%
L_{1}(\lambda )=(\cos \lambda )^{-1}$ is not defined outside of the interval 
$(-\frac{\pi }{2},\frac{\pi }{2})$. To end up this discussion, no mixing can
give increasing hazard rate for the NEF generated by the density $s_{1}(x)=%
\frac{1}{2\cosh \frac{\pi x}{2}}$.

\vspace{4mm}\textbf{\noindent }{\textbf{The hyperbolic case }}$\mathbf{%
\alpha =2}$: As we are going to see here, the situation is more favorable
when dealing with the following direct consequence of (\ref{HC}): 
\begin{equation*}
s_{2}(x)=\frac{x}{2\sinh \frac{\pi x}{2}}.
\end{equation*}%
Proposition 4 has shown that $b_{2}=-\log s_{2}$ is a convex function, so we
are in position to use Proposition 1. We have the explicit calculation 
\begin{equation*}
b_{2}^{\prime \prime }(x)=\frac{1}{x^{2}}-(\frac{\pi }{2})^{2}\frac{1}{\sinh
^{2}\frac{\pi x}{2}}.
\end{equation*}%
Since $|t|\leq |\sinh t|$, clearly $b_{2}^{\prime \prime }(x)>0$ and we have
a direct proof of the log concavity of $s_{2}$. Thus we have to study the
set of $(c,d)$'s such that the inequality 
\begin{equation}
cT(x)=c\frac{x\sinh \frac{\pi x}{2}}{\sqrt{\sinh ^{2}\frac{\pi x}{2}-(\frac{%
\pi x}{2})^{2}}}\leq \cosh (cx+d)  \label{Y}
\end{equation}%
holds for all real $x.$ For this we use the following lemma.

\begin{lemma}
\vspace{4mm}For all real $t$ we have 
\begin{equation*}
\frac{t\sinh t}{\sqrt{\sinh ^{2}t-t^{2}}}\leq \sqrt{3+t^{2}},
\end{equation*}%
where an equality occurs when $t=0.$
\end{lemma}

\begin{proof}
The proof follows from the inequality $\sinh ^{2}t-t^{2}-\frac{t^{4}}{3}\geq
0$, which is deduced from the expansion of 
\begin{equation*}
\sinh ^{2}t-t^{2}-\frac{t^{4}}{3}=\frac{1}{2}\cosh 2t-\frac{1}{2}-t^{2}-%
\frac{t^{4}}{3}=\sum_{n=3}^{\infty }\frac{2^{2n}t^{2n}}{(2n)!}\geq 0.
\end{equation*}
\end{proof}

\vspace{4mm}\noindent In terms of the function $T$, the lemma is equivalent
to stating that for all $x$ we have 
\begin{equation*}
T(x)\leq \frac{2}{\pi }\sqrt{3+(\frac{\pi x}{2})^{2}}
\end{equation*}%
Since it is difficult to find all $(c,d)$ such that $(\ref{P})$ holds we
shall content to study the set of $(c,d)$ such that $c\frac{2}{\pi }\sqrt{3+(%
\frac{\pi x}{2})^{2}}\leq \cosh (cx+d)$ holds for all $x$, or equivalently,
by introducing $k=2c/\pi $ and $u=k\frac{\pi x}{2}+d,$ to study the set of $%
(k,d)$ such that 
\begin{equation*}
\sqrt{3k^{2}+(u-d)^{2}}\leq \cosh u
\end{equation*}%
holds for all $u.$

\begin{lemma}
For $k>0$ and $d$ real $\sqrt{3k^{2}+(u-d)^{2}}\leq \cosh u$ holds for all $%
u $ if and only if $|d|\leq d_{0}=\sqrt{2}-\log (1+\sqrt{2})=0.532...$ and 
\begin{equation*}
3k^{2}\leq (2-\cosh ^{2}u_{d})\cosh ^{2}u_{d},
\end{equation*}%
where $u_{d}$ is the solution of the equation $\sinh 2u=2(u-d).$ In
particular, the inequality $\sqrt{3k^{2}+u^{2}}\leq \cosh u$ holds for all $%
u $ if and only if $k\leq \sqrt{\frac{2}{3}}.$
\end{lemma}

\begin{proof}
The inequality $\sqrt{3k^{2}+(u-d)^{2}}\leq \cosh u$ implies $|u-d|\leq
\cosh u$ for all $u.$ Now the minimum $d_{0}$ of $\cosh u-u$ is attained at $%
\log (1+\sqrt{2})$ which is the solution of the equation $\sinh u-1=0$, and
thus, $d_{0}=\sqrt{2}-\log (1+\sqrt{2}).$ Similarly, the minimum of $\cosh
u+u$ is attained at $-\log (1+\sqrt{2})$ and is $d_{0}.$ Since, $-\cosh
u-u\leq -d_{0}\leq d_{0}\leq \cosh u-u$, we get that $|u-d|\leq \cosh u$ for
all $u$ if and only if $|d|\leq d_{0}.$ Now, fix $d\in \lbrack -d_{0},d_{0}]$%
, then for finding all $k$ such that $3k^{2}+(u-d)^{2}\leq \cosh ^{2}u$, we
look for the (positive) minimum of $u\mapsto \cosh ^{2}u-(u-d)^{2}$, which
is attained at the point $u_{d}.$ Letting $d=0$ gives $u_{d}=0$ and also
entails the final result.
\end{proof}

\vspace{4mm}\noindent \textbf{Practical conclusion for 
$\alpha =2.$} The Laplace transform of $s_2$ is $1/\cos^2 \lambda$ for $\lambda\in (-\frac{\pi}{2},\frac{\pi}{2}).$ According to Lemma 7, we fix any $d$ such that $|d|\leq 0.532..$
and $c>0$ such that $c\leq \frac{\pi}{2}\sqrt{\frac{2}{3}}=1.282..$ We now choose an arbitrary number $\lambda$ such that
$$-\frac{\pi}{2}<\lambda_1=\lambda-c<\lambda+c=\lambda_2<\frac{\pi}{2}$$ and we use formula (\ref{PP}) for defining the mixing coefficient $p$ depending on $d$, $\lambda_1$ and $\lambda_2.$ With this choice, the density
$$\frac{x}{2\sinh \frac{\pi x}{2}}[pe^{-\lambda_1 x}\cos ^2\lambda_1+(1-p)e^{-\lambda_2 x}\cos ^2\lambda_2]$$
has an increasing hasard rate.

\section{The Ressel NEF}

Consider the density $s_{\alpha }$ on the positive real line defined by 
\begin{equation}
s_{\alpha }(x)=e^{-b_{\alpha }(x)}=\frac{\alpha \;x^{x+\alpha -1}e^{-x}}{%
\Gamma (x+\alpha +1)},\text{ }x>0,\text{ }\alpha >0,  \label{Ressel}
\end{equation}%
then $s_{\alpha }$ is called the Ressel or the Kendall-Ressel density with
parameter\ $\alpha $. It is infinitely divisible (in other terms the
Jorgensen set $J(s_1)$ is $[0,\infty)$) and 
\begin{equation*}
s_{\alpha}*s_{\alpha^{\prime}}=s_{\alpha+\alpha^{\prime}}.
\end{equation*}%
This density appears in various areas. For an M/G/1 queueing system with
arrival rate $\lambda $, it is the limiting distribution, as $\lambda
\rightarrow \infty $, of the length of the busy period $T(\alpha )-\alpha $
initiated by the virtual time quantity $\alpha >0 $ (c.f. Prabhu,1965, pages
73 and 237). In their characterization of the regression of the sample
variance on the sample mean, Fosam and Shanbhag (1997) showed that such a
regression is cubic on the sample mean for only for six distributions, of
which one is the Kendall-Ressel distribution. Kokonendji (2001) also
revealed this distribution on his investigation of first passage times on $0$
and $1$ of some L\'{e}vy processes for NEF's. Additional references can be
provided here regarding the Kendall-Ressel distribution like Pakes (1995)
formula (4.1), but one of the most detailed reference is Letac and Mora
(1990) who characterized all NEF's having cubic VF's, of which, of course,
the Ressel NEF is one of them.\ 

The Ressel NEF generated by (\ref{Ressel}) has a VF $(V,\Omega )=\left( 
\frac{\mu ^{2}}{\alpha}(1+\frac{\mu }{\alpha}),(0,\infty )\right) $.\ We are
interested in the values of $\alpha$ such that $b^{\prime\prime}_{%
\alpha}(x)\geq 0$ for all $%
x>0.$ One can consult Proposition 5.5 of Letac and Mora (1990) for checking
the puzzling formula $\int_{0}^{\infty }s_{\alpha}(x)dx=1$ and page 36 of
this reference for learning why this density can also be called the
Kendall-Ressel density.

\begin{proposition}
Let 
\begin{equation}
g(x)=\frac{\alpha-1}{x^{2}}+\frac{(2-\alpha)x-\alpha^{2}+\alpha}{x(x+\alpha)}%
\text{ and }h(x)=\frac{(\alpha^{2}-1)+(\alpha-\alpha^{2})x+(2-\alpha)x^{2}}{%
x^{2}(x+\alpha+1)},  \label{hr}
\end{equation}%
then for all $x>0$ we have $h(x)\leq
b^{\prime\prime}_{\alpha}(x)\leq g(x).$ Furthermore there exists a number $a\in (1.77,1.91)$ such
that $b^{\prime\prime}_{\alpha}(x)\geq 0$ for all $x>0$ if and only if $%
\alpha\in [ 1,a].$ In other terms the Glaser set $G(s_1)$ is $[
1,a].$
\end{proposition}

\begin{proof}
For $x>0$, we use the digamma function $\psi=\Gamma'/\Gamma$ and formula (\ref{DG}).
Using this notation we have 
\begin{eqnarray}
f(x)=b_{\alpha }^{\prime \prime }(x) &=&-\frac{\alpha -1}{x}+\frac{\alpha -1%
}{x^{2}}+\psi ^{\prime }(x+\alpha +1)  \notag \\
&=&-\frac{\alpha -1}{x}+\frac{\alpha -1}{x^{2}}+\sum_{n=2}^{\infty }\frac{1}{%
(n+x+\alpha -1)^{2}}.  \label{RA}
\end{eqnarray}%
Now observe that 
\begin{equation*} 
\frac{1}{1+x+\alpha }=\int_{2}^{\infty }\frac{dt}{(t+x+\alpha -1)^{2}}%
<\sum_{n=2}^{\infty }\frac{1}{(n+x+\alpha -1)^{2}}<\frac{1}{x+\alpha }%
=\int_{1}^{\infty }\frac{dt}{(t+x+\alpha -1)^{2}}.
\end{equation*}%
This gives the desired inequalities $h\leq f\leq g$. Clearly the function $f$
is positive if $\alpha =1$. If $\alpha <1,$ the function $f$ is equivalent
in a neighborhood of $x=0$ to $\frac{\alpha -1}{x^{2}}$ which tends to $%
-\infty $. It is obvious that for $\alpha \geq 2$, $g(x)$ becomes negative
ultimately (if $\alpha =2$ then $g(x)=(2-x)/x^{2}(2+x)$). Hence assume that $%
\alpha <2$ and let us study the sign of $g$. Since $g(x)=0$ if 
\begin{equation*}
(\alpha ^{2}-\alpha )+(-1+2\alpha -\alpha ^{2})x+(2-\alpha )x^{2}=0,
\end{equation*}%
then this equation has at least one solution if 
\begin{equation*}
D(\alpha )=(-1+2\alpha +\alpha ^{2})^{2}-4(2-\alpha )(\alpha ^{2}-\alpha
)=1+4\alpha -6\alpha ^{2}+\alpha ^{4}
\end{equation*}%
is nonnegative, which is the case for $\alpha >\alpha ^{\ast }=1.90321$. One
of the two possible solutions of $g(x)=0$ is 
\begin{equation*}
x_{1}=\frac{1-2\alpha +\alpha ^{2}+\sqrt{D(\alpha )}}{4-2\alpha }>\frac{%
(\alpha ^{\ast })^{2}-3}{4-2\alpha }>0,\quad \alpha ^{\ast }<\alpha <2.
\end{equation*}%
So $g(x_{1})=0$ and hence $g(x)<0$ for some $x>0$. On the other hand, we use
inequality $h\leq f$ and study the sign of $h$. Let $1<\alpha <2$, then if $%
D(\alpha )=(\alpha -1)(-8-4\alpha +3\alpha ^{2}+\alpha ^{3})<0$, $h(x)$ has
no zeros at all and hence $h(x)>0$ for all $x>0$. Now, $D(\alpha )=0$ if $%
\alpha \in \{-3.48929,-1.28917,1.0,1.77846\}$ and $D(\alpha )\rightarrow
\infty $ as $\alpha \rightarrow \infty $, hence $D(\alpha )<0$ (i.e., $f$ is
positive for all $x>0$)\ if $1<\alpha <1.77846$.
\end{proof}

\vspace{4mm}\noindent \textbf{Remarks.} For studying the log concavity of
the density $s_{\alpha }$ of the Ressel distribution one can be tempted to
imitate Proposition 4 and to wonder if for $\alpha >1$ the function 
\begin{equation*}
\frac{1}{s_{\alpha }(x)}=\frac{1}{\alpha }\frac{1}{x^{\alpha -1}}e^{x}\Gamma
(x+\alpha +1)\times x^{-x}
\end{equation*}%
defined on $(0,\infty )$ is the LT of a positive measure. Proposition 8 has
shown that is is impossible if $\alpha >a.$ This can be explained by the
fact that the factor $x\mapsto x^{-x}$ is the only factor in $1/s_{\alpha }$
which is not a LT (this observation follows from the fact that $x\mapsto
x^{x}$ is the LT of a stable law with parameter $1$, and the reciprocal of
the LT of a non Dirac measure cannot be a LT). In terms of the Glaser and
Karlin sets, the density $s_{1}$ is quite interesting. Proposition 8 has
shown that $G(s_{1})=[1,a]\subset K(s_{1})$ and a striking consequence is
that while the density $s_{1}$ is log concave the density $s_{2}=s_{1}\ast
s_{1}$ is not. This demonstrates the difference between the Glaser and the
Karlin sets. The additive semigroup generated by $[1,a]$ is $[1,a]\cup
\lbrack 2,\infty )\subset K(s_{1}).$ One can reasonably conjecture that $%
K(s_{1})=[1,\infty ).$

\vspace{4mm}\noindent From Proposition 8 it follows that if $\alpha\in
\lbrack 1,a]$, we are led to consider $T(x)=1/\sqrt{b^{\prime\prime}_{%
\alpha}(x)}$ and study the set of the pairs $(c,d)$  such that $%
cT(x)\leq \cosh (cx+d)$. Since $h(x)\leq
b^{\prime\prime}_{\alpha}(x)=-(s_{\alpha}^{\prime }/s_{\alpha})^{\prime }(x)$%
, a sufficient condition for these $(c,d)$ values is that $c/\sqrt{h(x)}\leq
\cosh (cx+d)$. Even this simplified inequality is still too complicated and
we shall content here to consider only the case $\alpha=1.$ For this case,
we search for the set of $(c,d)$'s such that for all $x>0$ we have $c\sqrt{%
x+2}\leq \cosh (cx+d).$ Therefore, the next proposition is devoted to the
case of the NEF generated by the probability density 
\begin{equation}  \label{Rone}
s_{1}(x)=\frac{x^{x}e^{-x}}{\Gamma (x+2)}\mathbf{1}_{(0,\infty )}(x)dx.
\end{equation}

\begin{proposition}
For $c>0$ and $d$ real consider the function 
\begin{equation*}
\varphi (x)=\frac{1}{c^{2}}\cosh ^{2}(cx+d)-x-2
\end{equation*}%
and define $x_{0}=\frac{1}{c}\left( \frac{1}{2}\log (c+\sqrt{c^{2}+1}%
)-d\right) .$ Then $\varphi (x)>0$ for all $x>0$ if and only if
\end{proposition}

\begin{enumerate}
\item either $x_0\leq 0$ and $\cosh d\geq \sqrt{2}c;$

\item or $x_{0}\geq 0$, $c\leq \sqrt{\frac{8}{7}}$ and $x_{0}\leq \frac{1}{%
2c^{2}}(1+\sqrt{c^{2}+1})-\frac{7}{4}.$
\end{enumerate}

\begin{proof}
Since $\varphi ^{\prime }(x)=\frac{1}{c}\sinh 2(cx+d)-1$, $\varphi ^{\prime
}(x_{0})\leq 0$ for $x\leq x_{0}$ and $\varphi ^{\prime }(x_{0})\geq 0$ for $%
x\geq x_{0}$, it follows that for $x_{0}\leq 0$, $\varphi (x)\geq 0$ for $%
x\geq 0$ if and only if $\varphi (0)\geq 0$ or if and only if $\cosh d\geq 
\sqrt{2}c$ (recall that $x_{0}\leq 0$ implies $d>0$). Similarly for $%
x_{0}\geq 0$, we have $\varphi (x)\geq 0$ for $x\geq 0$ if and only if $%
\varphi (x_{0})\geq 0$. Since $\frac{1}{c^{2}}\cosh ^{2}(cx_{0}+d)=\frac{1}{%
2c^{2}}(1+\sqrt{c^{2}+1}),$ the inequality $\varphi (x_{0})\geq 0$ is
equivalent to $x_{0}\leq \frac{1}{2c^{2}}(1+\sqrt{c^{2}+1})-\frac{7}{4}$,
which can be realized only if the right hand side is nonnegative, that is if 
$c\leq \sqrt{\frac{8}{7}}$.\noindent
\end{proof}

Now, having pairs $(c,d)$ to our disposal, we have to compute $p$ as
given by the formula (\ref{PP}) and we need for this the values of the
Laplace transform at points $\lambda _{1}=\lambda -c$ and $\lambda
_{2}=\lambda +c.$ For the Ressel distribution $s_{1}$, its LT 
\begin{equation}
L_{1}(\lambda )=\int_{0}^{\infty }e^{-\lambda x}s_{1}(x)dx  \label{LRone}
\end{equation}%
cannot be expressed explicitly. However, a numerical or graphical
calculation of $L_{1}(\lambda )$ for a given positive value of $\lambda $ is easily
done by means of the following proposition. Its statement is equivalent to
formula (11) in Fosam and Shanbhag (1997) which relies on Prabhu (1965, p.
73 and p. 237). We give an independent proof here for sake of completeness.

\begin{proposition}
Consider the bijection $f$ from $[1,\infty )$ to itself defined by $%
f(x)=x-\log x.$ Then for a given $\lambda >0$, the number $1/L_{1}(\lambda )$
defined by (\ref{LRone}) satisfies $f(1/L_{1}(\lambda ))=1+\lambda .$%
\noindent
\end{proposition}

\begin{proof}
Let $(Y(t))_{t\geq 0}$ be the L\'{e}vy process such that for $t,\lambda \geq
0$ we have $\mathbb{E}(e^{-\lambda Y(t)})=(1+\lambda )^{-t}$ (such a process
is usually called the gamma process). Define $T=\min \{t:t-Y(t)=1\}$. The
random variable $T-1$ has the Ressel distribution (\ref{Rone}) ( see also
Letac and Mora, 1990, p. 27 for this observation). Furthermore, Theorem 5.3
in Letac and Mora (1990) states that if $\mathbb{E}(e^{\theta
T})=e^{-q(\theta )}$ and $\mathbb{E}(e^{\theta (1-Y(1)})=e^{-r(\theta )}$
then $r(q(\theta ))=\theta $. Thus, since $L_{1}(\lambda )=\mathbb{E}%
(e^{-\lambda (T-1)})=e^{\lambda -q(-\lambda )}$ and since $e^{-r(\theta )}=%
\frac{e^{\theta }}{1+\theta }$, we can write $r(q(\theta ))=\log (1+q(\theta
))-q(\theta )=\theta $. This leads to 
\begin{equation*}
\log (1+q(-\lambda ))=q(-\lambda )-\lambda =-\log L_{1}(\lambda )
\end{equation*}%
The elimination of $q(-\lambda )$ between these two equalities gives $\frac{1%
}{L_{1}(\lambda )}+\log L_{1}(\lambda )=1+\lambda $. Since, for $\lambda
\geq 0$ we have $L_{1}(\lambda )\leq 1$, we get $f(1/L_{1}(\lambda
))=1+\lambda $.
\end{proof}

\vspace{4mm}\noindent \textbf{Practical conclusion for 
$\alpha =1.$} We fix a pair $(c,d)$ such that either condition 1 or condition 2 of Proposition 9 holds. We choose a number $\lambda$ such that $0<\lambda_1=\lambda-c<\lambda+c=\lambda_2.$ We compute numerically $L_1(\lambda_1)$ and $L_1(\lambda_2)$ with Proposition 10. The mixing coefficient $p$ is therefore determined by (\ref{PP}). The density on $(0,\infty)$
$$\frac{x^{x}e^{-x}}{\Gamma (x+2)}\left[p\frac{e^{-\lambda_1x}}{L_1(\lambda_1)}+(1-p)\frac{e^{-\lambda_2x}}{L_1(\lambda_2)}\right]$$
has an increasing hazard rate.

\section{The Kummer type 2 NEF}

Let $a,\lambda >0$ and $b$ real and consider the number 
\begin{equation*}
C(a,b,\lambda )=\int_{0}^{\infty }\frac{x^{a-1}}{(1+x)^{a+b}}e^{-\lambda
x}dx.
\end{equation*}%
As a function of $\lambda $, $C$  is proportional to what is sometimes
called in the literature the confluent hypergeometric function of the second
kind or a Whittaker function. If 
\begin{equation}
s(x)=\frac{1}{C(a,b,\lambda )}\frac{x^{a-1}}{(1+x)^{a+b}}e^{-\lambda x}%
\mathbf{1}_{(0,\infty )}(x),  \label{Ktwo}
\end{equation}%
then the probability $s(x)dx=K^{(2)}(a,b,\lambda )(dx)$ is called the Kummer
distribution of type 2 with parameters $(a,b,\lambda ).$ Needless to say,
that if $(a,b)$ are fixed the model $\{s(x)dx,\ \lambda >0\}$ is an NEF. If $%
b>0$, this model is generated by the beta distribution of type 2, i.e., by 
\begin{equation*}
\beta ^{(2)}(a,b)(dx)=\frac{1}{B(a,b)}\frac{x^{a-1}}{(1+x)^{a+b}}\mathbf{1}%
_{(0,\infty )}(x)dx.
\end{equation*}%
Kummer distributions have been studied by Ng and Kotz (1995). Statistical
aspects of Kummer distributions for waiting times and exceedance statistics
have been considered by Fitzgerald (2002). The Kummer distributions of type
one belong to NEF's generated by the ordinary beta distributions. Since they
are concentrated on the bounded set $(0,1)$, they are not relevant for our
study. Accordingly, we study the NEF generated the Kummer distribution of
type 2. Its VF cannot be expressed explicitly. However, the important fact
about such an NEF is the formula (\ref{Dyson}) below which gives the LT of $%
\beta ^{(2)}(a,b)$ in terms of the confluent hypergeometric function defined
for real $a$ and $b$ such that $b$ is not in the set $-\mathbb{N}$. This LT
is then given in terms of the entire function 
\begin{equation}
_{1}F_{1}(a\;;b\;;\lambda )=\sum_{n=0}^{\infty }\frac{(a)_{n}\lambda ^{n}}{%
n!(b)_{n}}.  \label{KM}
\end{equation}%
Here, $(a)_{0}=1$ and $(a)_{n+1}=(a+n)(a)_{n}$. This formula states that if $%
a>0$ and $b$ is not in the set $\mathbb{Z}$ of relative integers, we have 
\begin{equation}
C(a,b,\lambda )=\frac{\Gamma (b)\Gamma (a)}{\Gamma (a+b)}\
_{1}F_{1}(a\;;1-b\;;\lambda )+\Gamma (-b)\lambda ^{b}\
_{1}F_{1}(a+b\;;1+b\;;\lambda ).  \label{Dyson}
\end{equation}%
In (\ref{Dyson}), the mapping $z\mapsto 1/\Gamma (z)$ is an entire function
which coincides with the ordinary $1/\Gamma (z),z>0,$ and (\ref{Dyson}) can
be extended to the case where $b\in \mathbb{Z}$ by a limiting process. The
identity (\ref{Dyson}) is by no means elementary and its proof by the Barnes
formula can be found for instance in Slater (1960), formula 3.1.19. A
probabilistic proof would be desirable.

To exemplify the use of (\ref{Dyson}), observe that if $a,b,\lambda >0$
and if $X\sim \gamma(b,\lambda ),\ Y\sim K^{(2)}(a,b,\lambda )$ and $Z\sim
\gamma(a+b,\lambda )$ are independent, then $X+Y$ and $Z/(1+Y)$ have the
same distribution $K^{(2)}(a+b,-b,\lambda ).$ In order to prove this, for
suitable $t^{\prime }s$ just consider the LT $\mathbb{E}(e^{-t(X+Y)})$ and
the Mellin transform $\mathbb{E}(Z^{t}/(1+Y)^{t}),$ and then use (\ref{Dyson}%
).

\begin{proposition}
If $s$ is defined by (\ref{Ktwo}) then $-(s^{\prime }/s)^{\prime }>0$ for
all $x>0$ if and only if $1\leq a$ and $b\leq -1$ with $a-b-2\neq 0.$
\end{proposition}

\begin{proof}
If $A=a-1$ and $B=-b-1$, we get 
\begin{equation}
-(\frac{s^{\prime }}{s})^{\prime }(x)=\frac{A+2Ax+Bx^{2}}{x^{2}(1+x)^{2}}.
\label{KP}
\end{equation}%
Trivially, $-(s^{\prime }/s)^{\prime }>0$ if $A\geq 0$ and $B\geq 0$ with $%
A+B\neq 0.$ Conversely, if $A+2Ax+Bx^{2}>0$ for all $x>0$, then letting $%
x\rightarrow \infty $ shows that $B\geq 0$. Also, letting $x\rightarrow 0$
shows that $A\geq 0,$ while $A=B=0$ would imply that $-(s^{\prime
}/s)^{\prime }=0.$\vspace{4mm}Therefore if $A=a-1\geq 0$ and $B=-b-1$ with $%
AB\neq 0$, we are allowed, by using (\ref{KP}), to consider for $x>0$, 
\begin{equation*}
T(x)=1/\sqrt{-(\frac{s^{\prime }}{s})^{\prime }(x)}=\frac{x+x^{2}}{\sqrt{%
A+2Ax+Bx^{2}}}.
\end{equation*}%
We then have to investigate which numbers $c>0$ and $d$ real are such that
for all $x>0,$ we have 
\begin{equation*}
c\frac{x+x^{2}}{\sqrt{A+2Ax+Bx^{2}}}\leq \cosh (cx+d)
\end{equation*}%
For simplicity, we are going to treat only the particular case $A=0$, and
therefore to study the NEF 
\begin{equation*}
\frac{1}{C(1,-B-1,\lambda )}(1+x)^{B}e^{-\lambda x},
\end{equation*}%
where $B$ is a fixed positive constant. For this particular case we are
looking for the values of $(c,d)$ with $c>0$ such that for all $x>0$ we have 
$\frac{\sqrt{B}}{c}\cosh (cx+d)-x-1\geq 0.$
\end{proof}

\begin{proposition}
For $B,c>0$ and $d$ real, consider the function defined on $\mathbb{R}$ by 
\begin{equation*}
\varphi (x)=\frac{\sqrt{B}}{c}\cosh (cx+d)-x-1
\end{equation*}%
and define $x_{0}=\frac{1}{c}\left( \log \frac{1+\sqrt{B+1}}{\sqrt{B}}%
-d\right) $. Then $\varphi (x)\geq 0$ for all $x\geq 0$ if and only
\end{proposition}

\begin{enumerate}
\item either $x_0\leq 0,$ $\sqrt{B}\leq c$ and $\cosh d\geq \frac{c}{\sqrt{B}%
};$

\item or $x_{0}\geq 0$, $c\leq \sqrt{B+1}$ and $x_{0}\leq 1-\frac{\sqrt{B+1}%
}{c}.$
\end{enumerate}

\begin{proof}
We study the function $\varphi $ in an elementary way: We get that $\varphi
^{\prime }(x)=\sqrt{B}\sinh (cx+d)-1$ satisfies $\varphi ^{\prime }(x)\leq 0$
for $x\leq x_{0}$ and $\varphi ^{\prime }(x)\geq 0$ for $x\geq x_{0}$. If $%
x_{0}\leq 0$ then $\varphi (x)\geq 0$ for all $x\geq 0$ if and only $\varphi
(0)=\frac{\sqrt{B}}{c}\cosh d-1\geq 0$ and this proves part 1. If $x_{0}\geq
0$, then $\varphi (x)\geq 0$ for all $x\geq 0$ if and only 
\begin{equation*}
\varphi (x_{0})=\frac{\sqrt{B+1}}{c}-1-x_{0}\geq 0
\end{equation*}%
which proves part 2.
\end{proof}

\vspace{4mm}\noindent Here, again, in order to apply the results of this
section to formula (\ref{PP}) we have to compute the values of the LT $%
C(1,-B-1,\lambda ),$ which can also be seen as a truncated gamma function.
If $B$ is an integer, $C$ is easily computed by the binomial formula 
\begin{equation*}
C(1,-B-1,\lambda )=\int_{0}^{\infty }(1+x)^{B}e^{-\lambda x}dx=\frac{B!}{%
\lambda ^{B+1}}\sum_{n=0}^{B}\frac{\lambda ^{n}}{n!}.
\end{equation*}%
If $B>0$ is not an integer, (\ref{Dyson}) gives 
\begin{eqnarray*}
C(1,-B-1,\lambda ) &=&e^{\lambda }\int_{1}^{\infty }x^{B}e^{-\lambda x}dx=%
\frac{\Gamma (B+1)}{\lambda ^{B+1}}e^{\lambda }-\frac{1}{B+1}\
_{1}F_{1}(1;2+B;\lambda ) \\
&=&\frac{\Gamma (B+1)}{\lambda ^{B+1}}\left[ e^{\lambda }-\sum_{n=0}^{\infty
}\frac{\lambda ^{B+n+1}}{\Gamma (B+n+2)}\right] ,
\end{eqnarray*}%
but then we have to rely on numerical analysis for computing the
corresponding values of the confluent hypergeometric function 
\begin{equation*}
_{1}F_{1}(1;B+2;\lambda )=1+\sum_{n=1}^{\infty }\frac{\lambda ^{n}}{%
(B+2)\ldots (B+n+1)} 
\end{equation*}%
and use (\ref{Dyson}). A good reference for such numerical consideration aspects can be found in
Abad and Sesma (1995).

\vspace{4mm}\noindent \textbf{Practical conclusion for 
$a=1, b=-1-B.$} We fix a pair $(c,d)$ such that either condition 1 or condition 2 of Proposition 12 holds. We choose a number $\lambda$ such that $0<\lambda_1=\lambda-c<\lambda+c=\lambda_2.$ We compute numerically $C(1;-1-B,\lambda_1)$ and $C(1;-1-B,\lambda_2)$. The mixing coefficient $p$ is therefore determined by (\ref{PP}). The density on $(0,\infty)$
$$(1+x)^B\left[p\frac{e^{-\lambda_1x}}{C(1;-1-B,\lambda_1)}+(1-p)\frac{e^{-\lambda_2x}}{C(1;-1-B,\lambda_2)}\right]$$
has an increasing hazard rate.

\end{document}